%% file: main_ArXiv.tex
\documentclass{article}
\usepackage{arxiv}
\include{header}

\mathtoolsset{showonlyrefs}

\floatname{algorithm}{Алгоритм}
\newtheorem{theorem}{Теорема}
\theoremstyle{definition}
\newtheorem{definition}{Определение}
\newtheorem{assumption}{Предположение}
\newtheorem{remark}{Замечание}
\renewcommand{\leq}{\leqslant}
\renewcommand{\geq}{\geqslant}
\DeclarePairedDelimiter{\rbr}{(}{)}

\DeclarePairedDelimiter{\cbr}{\{}{\}}
\DeclarePairedDelimiter{\abr}{\langle}{\rangle}
\DeclarePairedDelimiter{\norm}{\lVert}{\rVert}
\DeclarePairedDelimiter{\ceil}{\lceil}{\rceil}

\begin{document}

\title{Градиентные методы для минимизационных задач с условием Поляка--Лоясиевича: относительная погрешность градиента и адаптивный подбор параметров}

\author{
    Пучинин Сергей Максимович \\
    МФТИ \\
    Москва, Россия \\
    \texttt{puchinin.sm@phystech.edu} \\
    \And
    Стонякин Федор Сергеевич \\
    МФТИ \\
    Москва, Россия \\
    КФУ им. В.\,И.\,Вернадского \\
    Симферополь, Россия \\
    \texttt{fedyor@mail.ru}
}

{\raggedright УДК~519.85}\bigskip

\date{}

\maketitle

\begin{abstract}
\end{abstract}

\section{Введение}

С ростом набора прикладных задач, которые сводятся к оптимизационным задачам большой или даже огромной размерности (некоторые из таких приложений возникают в машинном обучении, глубоком обучении, оптимальном управлении, обработке сигналов, статистике и т.д.), исследование различных методов первого порядка привлекает большое внимание научного сообщества \cite{beck2017first}. Градиентные методы можно рассматривать как одно из ключевых направлений развития численных методов оптимизации в настоящее время. Преимущественно это связано с малозатратностью их итераций с точки зрения используемой памяти, а также с независимостью оценок скорости сходимости от параметров размерности пространства.

Для задачи минимизации гладкой функции $f$ хорошо известно, что если $f$ сильно выпукла, то метод градиентного спуска имеет глобальную линейную скорость сходимости \cite{nesterov2004introductory}, причём её оценки не зависят от размерности пространства. Однако многие фундаментальные задачи машинного обучения, такие как регрессия методом наименьших квадратов или логистическая регрессия приводят к задачам минимизации функций, которые не сильно выпуклы и даже не выпуклы. Это привело к исследованию приемлемых для возникающих в такого типа приложениях аналогов свойств выпуклости и сильной выпуклости для целевой функции оптимизационной задачи. Одно из наиболее известных таких свойств~--- условие градиентного доминирования Поляка--Лоясиевича ~\cite{polyak1963gradient}. Хорошо известно ~\cite{polyak1963gradient}, что этого условия достаточно, чтобы показать глобальную линейную скорость сходимости градиентного спуска для достаточно гладких задач без предположений о выпуклости. Недавно в \cite{Yue2022low} показана оптимальность этого результата.

Как правило, для методов первого порядка делается предположение о доступности в произвольной точке допустимого множества задачи точного оракула первого порядка, т.е. оракул должен выдавать в каждой запрашиваемой точке точные значения целевой функции и её градиента. Но, к сожалению, во многих приложениях нет возможности получить точную информацию о градиенте и/или целевой функции на каждой итерации метода. Это привело исследователей к изучению поведения методов первого порядка, которые могут работать с неточным оракулом. В работе \cite{devolder2013first} (которую можно считать фундаментальной в этом направлении) авторы вводят понятие неточного оракула первого порядка, которое естественно возникает во многих ситуациях. При этом вопросы исследования влияния погрешностей доступной информации на гарантии сходимости численных методов, по-видимому, впервые были поставлены в книге \cite{polyak1983introduction}.

Часто при анализе сходимости метода градиентного спуска подразумевается постоянная величина размера шага, которая зависит от константы Липшица градиента целевой функции (константы гладкости). Однако во многих прикладных задачах эту константу трудно оценить. Например, известная функция Розенброка и ее многомерные аналоги (например, функция Нестерова--Скокова) имеют только локально липшицево-непрерывный градиент. Для того чтобы преодолеть трудности с определением значения константы Липшица градиента, были разработаны различные методы, одним из которых является градиентный спуск с адаптивной политикой размера шага.

Недавно в \cite{polyak2022stopping} для рассматриваемого в настоящей работе класса задач (c достаточно гладкой целевой функцией, удовлетворяющей условию Поляка--Лоясиевича) предложен неадаптивный и адаптивный градиентные методы, использующие понятие неточного градиента. В \cite{polyak2022stopping} проанализированы предложенные алгоритмы и влияние помех в градиенте на скорость сходимости. Однако в их работе адаптивность имеет место лишь в отношении константы Липшица градиента; по-прежнему необходимо точно знать оценку величины абсолютной погрешности градиента. В \cite{kuruzov2023gradient} предложен и проанализирован адаптивный алгоритм, который предполагает настройку не только константы гладкости функции $L$, но и величины абсолютной погрешности градиента $\Delta$. Но в обеих работах \cite{kuruzov2023gradient,polyak2022stopping} рассматривался лишь случай доступности градиента с абсолютной погрешностью. Пример выбора размера шага неадаптивного градиентного метода в случае относительной погрешности задания градиента в случае известного значения $L$ приведён в параграфе 1 пособия \cite{gasnikov2021mccme}.

В настоящей работе продолжаются исследования по адаптивным методам градиентного типа \cite{polyak2022stopping, kuruzov2023gradient} для гладких задач с условием Поляка--Лоясиевича в случае доступности методу в произвольной текущей точке информации о градиенте с относительной погрешностью. Предложены два адаптивных алгоритма (различающихся набором настраиваемых параметров) для задач с целевыми функциями, удовлетворяющими условию Поляка--Лоясиевича, при наличии относительной неточности задания градиента, с подробным анализом их скоростей сходимости и поведения образуемых ими траекторий. Если в первом из этих алгоритмов (алгоритме~\ref{alg:1}) адаптивность реализована только для параметра $L$, то во втором (алгоритме~\ref{alg:2}) имеется адаптивность по обоим параметрам: константе Липшица градиента целевой функции $L$ и величине относительной погрешности градиента $\alpha$. Таким образом, второй из предложенных алгоритмов полностью адаптивен.

\section{Постановка задачи и основные понятия}\label{sec:2}

Рассматривается минимизационная задача (в общем случае невыпуклая)
\begin{equation}\label{eq:problem}
    \min\limits_{x \in \mathbb{R}^n} f(x),
\end{equation}
где целевая функция является $L$-гладкой и удовлетворяет условию Поляка--Лоясиевича.

\begin{definition}
Дифференцируемая функция $f: \mathbb{R}^n \to \mathbb{R}$ называется $L$-гладкой относительно нормы $\lVert \cdot \rVert$, если для некоторой константы $L > 0$ выполнено
    \begin{equation}\label{eq:L}
        \lVert \nabla f(y) - \nabla f(x) \rVert \leq L \lVert y - x \rVert.
    \end{equation}
\end{definition}

\noindent Норма в данной работе всюду подразумевается евклидовой. Как хорошо известно, условие $L$-гладкости \eqref{eq:L} влечёт следующее неравенство
\begin{equation}\label{eq:L_alt}
    f(y) \leq f(x) + \abr*{\nabla f(x), y - x} + \frac{L}{2} \lVert y - x \rVert^2.
\end{equation}

\begin{definition}
    $L$-гладкая функция $f$ удовлетворяет условию Поляка--Лоясиевича (или, для краткости, PL-условию), если выполнено следующее неравенство
    \begin{equation}\label{eq:PL}
        f(x) - f^* \leq \frac{1}{2 \mu} \lVert \nabla f(x) \rVert^2 \quad \forall x \in \mathbb{R} ^n,
    \end{equation}
    где $\mu > 0$~--- некоторая константа, а $f^* \coloneqq f(x^*)$, где $x^* \in X^*$ ($X^*$~--- множество точных решений рассматриваемой минимизационной задачи).
\end{definition}
\noindent В литературе условие \eqref{eq:PL} также иногда называется условием градиентного доминирования.

В данной работе предлагается рассмотреть проблему поведения методов градиентного типа для задачи \eqref{eq:problem} в случае относительных помех информации о градиенте. Недавно в \cite{kuruzov2023gradient,polyak2022stopping} были детально исследованы адаптивные методы градиентного типа при наличии абсолютных помех. Мы же рассматриваем ситуацию, когда градиент известен с относительной погрешностью, то есть
\begin{equation}\label{eq:rel_err}
    \lVert \Tilde{\nabla} f(x) - \nabla f(x) \rVert \leq \alpha \lVert \nabla f(x) \rVert,
\end{equation}
где $\Tilde{\nabla} f(x)$~--- неточный градиент (доступен методу), а $\alpha \in [0, 0.5)$~--- некоторая константа, отвечающая за величину помех в градиенте. Данное условие на неточный градиент было введено и изучено в работах \cite{carter1991global,polyak1983introduction}. Обычно рассматриваются $\alpha$ из большего полуинтервала $[0,1)$, однако по причинам, описанным в разделе~\ref{sec:3}, нами рассматривается усеченный полуинтервал $[0,0.5)$. Отметим, что относительные помехи в градиенте могут возникать из-за приближенного его вычисления, например, в безградиентной оптимизации (см. \cite{berahas2022theoretical}).

Как легко видеть из \eqref{eq:rel_err},
\begin{equation}\label{eq:norm_ineqs}
    (1 - \alpha) \lVert \nabla f(x) \rVert \leq \lVert \Tilde{\nabla} f(x) \rVert \leq (1 + \alpha) \lVert \nabla f(x) \rVert,
\end{equation}
откуда можно получить условие типа Поляка--Лоясиевича \eqref{eq:PL} для неточного градиента
\begin{equation}\label{eq:PL_rel}
    f(x) - f^* \leq \frac{1}{2 \mu (1 - \alpha)^2} \lVert \Tilde{\nabla} f(x) \rVert^2.
\end{equation}

\section{Градиентный спуск с адаптивной политикой размера шага}\label{sec:3}

Задачу \eqref{eq:problem} предлагается решать методом градиентного спуска в виде
\begin{equation}\label{eq:method}
    x^{k+1} = x^k - h_k \Tilde{\nabla} f(x^k), \quad h_k > 0,
\end{equation}
где размер шага $h_k$ может зависеть от $L$ и $\alpha$, а также, в случае когда точное значение одного из них или их обоих неизвестно, от $L_{k+1}$ и $\alpha_{k+1}$.

Если при реализации метода \eqref{eq:method} градиент доступен с известной относительной погрешностью ${\alpha \in [0, 1)}$ и известен параметр $L > 0$, то выбор постоянного размера шага
\begin{equation}
    h_k = h = \frac{1}{L} \frac{(1 - \alpha)}{(1 + \alpha)^2}
\end{equation}
приводит к следующему результату (см. параграф~1 из пособия \cite{gasnikov2021mccme} и имеющиеся там ссылки):
\begin{equation}\label{eq:funcPLconst}
    f(x^{k+1}) - f(x^k) \leq -\frac{1}{2L} \frac{(1 - \alpha)^2}{(1 + \alpha)^2} \norm*{\nabla f(x^k)}^2.
\end{equation}
Пользуясь условием Поляка--Лоясиевича \eqref{eq:PL}, имеем следующую оценку на скорость сходимости по функции
\begin{equation}\label{eq:known_rate}
    f(x^N) - f^* \leq \rbr*{1 - \frac{\mu}{L} \frac{(1 - \alpha)^2}{(1 + \alpha)^2}}^N \rbr*{f(x^0) - f^*}.
\end{equation}

Однако цель данной работы~--- предложить аналог алгоритмам 1 и 2 из \cite{kuruzov2023gradient} для случая относительной погрешности градиента с адаптивной настройкой параметров $L$ и $\alpha$. По сути, необходимо предложить критерий выхода из итерации, содержащий норму неточного градиента $\lVert \Tilde{\nabla} f(x^k) \rVert$. Это обстоятельство затрудняет использование подхода с оценками типа \eqref{eq:funcPLconst} для нормы точного градиента. Поэтому рассмотрим альтернативный вариант выбора размера шага для метода \eqref{eq:method} с относительной погрешностью задания градиента, который позволит получить приемлемый аналог оценки  \eqref{eq:funcPLconst} для квадрата нормы неточного градиента.
Используя \eqref{eq:rel_err}, \eqref{eq:norm_ineqs} выпишем следующий аналог неравенства \eqref{eq:L_alt}:
\begin{equation}
    \begin{aligned}
        f(y)
        &\leq f(x) + \langle \nabla f(x), y - x \rangle + \frac{L}{2} \lVert y - x \rVert^2 \\
        &= f(x) + \langle \Tilde{\nabla} f(x), y - x \rangle + \frac{L}{2} \lVert y - x \rVert^2 + \langle \nabla f(x) - \Tilde{\nabla} f(x), y - x \rangle \\
        &\leq f(x) + \langle \Tilde{\nabla} f(x), y - x \rangle + \frac{L}{2} \lVert y - x \rVert^2 + \lVert \nabla f(x) - \Tilde{\nabla} f(x) \rVert \lVert y - x \rVert \\
        &\leq f(x) + \langle \Tilde{\nabla} f(x), y - x \rangle + \frac{L}{2} \lVert y - x \rVert^2 + \alpha \lVert \nabla f(x) \rVert \lVert y - x \rVert \\
        &\leq f(x) + \langle \Tilde{\nabla} f(x), y - x \rangle + \frac{L}{2} \lVert y - x \rVert^2 + \frac{\alpha}{1 - \alpha} \lVert \Tilde{\nabla} f(x) \rVert \lVert y - x \rVert,
    \end{aligned}
\end{equation}
т.е.
\begin{multline}\label{eq:cond}
    f(x^{k+1}) \leq f(x^k) + \abr*{\Tilde{\nabla} f(x^k), x^{k+1} - x^k} + \frac{L}{2} \norm*{x^{k+1} - x^k}^2 +\\+ \frac{\alpha}{1 - \alpha} \norm*{\Tilde{\nabla} f(x^k)} \norm*{x^{k+1} - x^k}.
\end{multline}
На базе неравенства \eqref{eq:cond} ниже будут предложены критерии выхода из итерации в алгоритмах \ref{alg:1} и \ref{alg:2}. Это неравенство гарантирует отсутствие зацикливания (выход из итерации в некоторый момент) соответствующих алгоритмов для $L$-гладких задач.

\subsection{Алгоритм с адаптивной настройкой параметра гладкости}\label{subsec:3.1}

Ниже приведен адаптивный алгоритм \ref{alg:1} с настройкой параметра гладкости $L$ при известном $\alpha \in [0; 0.5)$.
\begin{algorithm}
    \caption{Градиентный спуск с адаптивной настройкой $L$}
    \label{alg:1}
    \begin{algorithmic}
        \State
        \textbf{Вход:} $x^0$, $L_{min} \geq \mu > 0$, $L_0 \geq L_{min}$, $\alpha \in [0, 0.5)$.
        \begin{enumerate}
            \renewcommand{\labelenumi}{\textbf{\arabic{enumi}:}}
            \item
                $k = 0$.
            \item
                $L_{k+1} = \max\cbr*{\frac{L_k}{2}, L_{min}}$.
            \item
                \begin{equation}\label{eq:alg1_step}
                    x^{k+1} = x^k - \frac{1}{L_{k+1}} \frac{1 - 2\alpha}{1 - \alpha} \Tilde{\nabla} f(x^k)
                \end{equation}
            \item Если
                \begin{multline}\label{eq:alg1_cond}
                    f(x^{k+1}) \leq f(x^k) + \abr*{\Tilde{\nabla} f(x^k), x^{k+1} - x^k} + \frac{L_{k+1}}{2} \norm*{x^{k+1} - x^k}^2 +\\+ \frac{\alpha}{1 - \alpha} \norm*{\Tilde{\nabla} f(x^k)} \norm*{x^{k+1} - x^k},
                \end{multline}
                то переходим на шаг 5. Иначе $L_{k+1} \coloneqq 2 L_{k+1}$ и возвращаемся на шаг 3.
            \item
                Если правило остановки не выполнено, то $k \coloneqq k+1$ и переходим на шаг 2.
        \end{enumerate}
        \textbf{Выход:} $x^k$.
    \end{algorithmic}
\end{algorithm}
Проведем анализ скорости сходимости и траектории для предложенного алгоритма. Из \eqref{eq:method} имеем
\begin{equation}
    \norm*{x^{k+1} - x^k} = h_k \norm*{\Tilde{\nabla} f(x^k)}
\end{equation}
и
\begin{equation}
    \abr*{\Tilde{\nabla} f(x^k), x^{k+1} - x^k} = -h_k \norm*{\Tilde{\nabla} f(x^k)}^2.
\end{equation}
Объединяя полученные выражения с \eqref{eq:alg1_cond}, получаем
\begin{equation}\label{eq:cond_for_alpha}
    f(x^{k+1}) - f(x^k) \leq \rbr*{- h_k + \frac{L_{k+1} h_k^2}{2} + \frac{\alpha h_k}{1 - \alpha}} \norm*{\Tilde{\nabla} f(x^k)}^2.
\end{equation}
Так как $h_k$ не может быть отрицательным, минимум выражения, стоящего в скобках, достигается при
\begin{equation}\label{eq:h_k_2}
    h_k = \max\cbr*{0, \frac{1}{L_{k+1}} \frac{1 - 2\alpha}{1 - \alpha}}.
\end{equation}
Если $\alpha \in [0.5, 1)$, то минимум выражения, стоящего в скобках в \eqref{eq:cond_for_alpha}, равен $0$ при $h_k = 0$, что говорит об отсутствии гарантий сходимости метода \eqref{eq:method} вне зависимости от выбора $h_k$. Поэтому в \eqref{eq:rel_err}, как и в алгоритмах \ref{alg:1} и \ref{alg:2}, рассматриваются лишь $\alpha \in [0, 0.5)$.

При $\alpha \in [0, 0.5)$ и задании размера шага согласно \eqref{eq:alg1_step} имеем
\begin{equation}\label{eq:funcdiff}
    f(x^{k+1}) - f(x^k) \leq - \frac{1}{2L_{k+1}} \frac{(1 - 2\alpha)^2}{(1 - \alpha)^2} \norm*{\Tilde{\nabla} f(x^k)}^2.
\end{equation}
Объединив последнее неравенство с \eqref{eq:PL_rel}, получаем
\begin{equation}\label{eq:prefinal}
    f(x^{k+1}) - f^* \leq \rbr*{1 - \frac{\mu}{L_{k+1}} (1 - 2\alpha)^2} \rbr*{f(x^k) - f^*} = \rbr*{1 - \frac{\mu}{L_{k+1}} \xi} \rbr*{f(x^k) - f^*},
\end{equation}
где
\begin{equation}\label{eq:xi}
    \xi \coloneqq (1 - 2\alpha)^2.
\end{equation}
Можно считать, что $L_0 \leq L$, поэтому $\max_{k<N} L_{k+1} \leq 2L$. Следовательно,
\begin{equation}\label{eq:final_1}
    f(x^N) - f^* \leq \prod\limits_{k = 0}^{N - 1} \rbr*{1 - \frac{\mu}{L_{k+1}} \xi} \rbr*{f(x^0) - f^*} \leq \rbr*{1 - \frac{\mu}{L_{max}} \xi}^N \rbr*{f(x^0) - f^*},
\end{equation}
где
\begin{equation}\label{eq:L_max_1}
    L_{max} \coloneqq 2L.
\end{equation}
Получаем сходимость по функции со скоростью геометрической прогрессии. Это означает, что для достижения $\frac{\epsilon}{\mu}$-точности по функции требуется не более
\begin{equation}
    N_* = \ceil*{\frac{L_{max}}{\mu\xi} \log\rbr*{\frac{\mu \rbr*{f(x^0) - f^*}}{\varepsilon}}}
\end{equation}
итераций. C другой стороны, в силу PL-условия \eqref{eq:PL_rel}, эта же точность по функции может быть гарантирована при выполнении следующего правила остановки
\begin{equation}
    \norm*{\Tilde{\nabla} f(x^k)}^2 \leq 2 \epsilon (1 - \alpha)^2.
\end{equation}

Далее, получим ограничение на траекторию $\{x^k\}_{k=0}^N$, образуемую алгоритмом \ref{alg:1}.
\begin{equation}
    \frac{1}{2L_{k+1}} \frac{(1 - 2\alpha)^2}{(1 - \alpha)^2} \norm*{\Tilde{\nabla} f(x^k)}^2 \leq f(x^k) - f(x^{k+1})
\end{equation}
и, следовательно, согласно \eqref{eq:alg1_step},
\begin{equation}
    \begin{gathered}
        \norm*{x^{k+1} - x^k}^2 = \rbr*{\frac{1}{L_{k+1}} \frac{1 - 2\alpha}{1 - \alpha}}^2 \norm*{\Tilde{\nabla} f(x^k)}^2 \leq \frac{2}{L_{k+1}} \rbr*{f(x^k) - f(x^{k+1})} \leq\\
        \leq \frac{2}{L_{k+1}} \rbr*{f(x^k) - f^*} \leq \frac{2}{L_{min}} \rbr*{1 - \frac{\mu}{L_{max}} \xi}^k \rbr*{f(x^0) - f^*}.
    \end{gathered}
\end{equation}
Наконец,
\begin{equation}\label{eq:traj}
    \begin{gathered}
        \norm*{x^N - x^0} \leq \sum\limits_{k = 0}^{N - 1} \norm*{x^{k+1} - x^k} \leq \sqrt{\frac{2}{L_{min}} \rbr*{f(x^0) - f^*}} \sum\limits_{k = 0}^{N - 1} \rbr*{1 - \frac{\mu}{L_{max}} \xi}^{k/2} =\\
        = \sqrt{\frac{2}{L_{min}} \rbr*{f(x^0) - f^*}} \frac{1 - \rbr*{1 - \frac{\mu}{L_{max}} \xi}^{N/2}}{1 - \rbr*{1 - \frac{\mu}{L_{max}} \xi}^{1/2}} \leq \frac{2L_{max}}{\mu\xi} \sqrt{\frac{2}{L_{min}} \rbr*{f(x^0) - f^*}}.
    \end{gathered}
\end{equation}

Таким образом, верна следующая
\begin{theorem}\label{th:1}
     Пусть функция $f$ удовлетворяет условию Липшица градиента \eqref{eq:L} и условию Поляка--Лоясиевича \eqref{eq:PL}. И пусть при работе алгоритма~\ref{alg:1} выполнилось одно из следующих двух условий:
     \begin{enumerate}
         \item алгоритм~\ref{alg:1} проработал $N_*$ итераций, где
             \begin{equation}\label{eq:N_1}
                 N_* = \ceil*{\frac{L_{max}}{\mu\xi} \log\rbr*{\frac{\mu \rbr*{f(x^0) - f^*}}{\varepsilon}}};
             \end{equation}
         \item для некоторого $N < N_*$ выполнено правило остановки в виде
            \begin{equation}
                \norm*{\Tilde{\nabla} f(x^N)}^2 \leq 2 \epsilon (1 - \alpha)^2,
            \end{equation}
     \end{enumerate}
     где $L_{max}$ задается формулой \eqref{eq:L_max_1}. Тогда верны следующие оценки на выходную точку $\hat{x}$ алгоритма~\ref{alg:1} ($\hat{x} = x^{N_*}$ или $\hat{x} = x^{N}$ соответственно):
     \begin{equation}
         f(\hat{x}) - f^* \leq \frac{\epsilon}{\mu}
     \end{equation}
     и
     \begin{equation}\label{eq:traj_1}
         \norm*{\hat{x} - x^0} \leq \frac{2 L_{max}}{\mu\xi} \sqrt{\frac{2}{L_{min}} \rbr*{f(x^0) - f^*}}.
     \end{equation}
\end{theorem}

\begin{remark}\label{Puch_rem1}
Оценка \eqref{eq:traj} не зависит ни от $N$, ни от невязки по функции на $N$-ом шаге. Поэтому для $L$-гладких задач траектория рассмотренного алгоритма лежит в некотором шаре фиксированного радиуса с центром в $x^0$. Данная особенность позволяет в какой-то мере сохранять полезные свойства начальной точки $x^0$. Также она позволяет расширить область применимости рассматриваемого подхода на другие классы функций. Например, на класс функций, удовлетворяющих PL-условию лишь локально на некотором множестве. Это продемонстрировано ниже, в вычислительных экспериментах в разделе~\ref{sec:5}.
\end{remark}

\begin{remark}
    Во всех приведённых в настоящем подразделе статьи результатах $L_{max}$, оценивающее максимально возможное значение для $L_{k+1}$, может быть заменено на $\max_{k < N} L_{k+1}$.
\end{remark}

\begin{remark}
    Алгоритм \ref{alg:1} применим и для случая точного градиента при $\alpha = 0$. В этом случае просто $\xi = 1$.
\end{remark}

\subsection{Алгоритм с адаптивной настройкой параметра гладкости и величины относительной погрешности градиента}\label{subsec:3.2}

Для реализации адаптивной настройки не только параметра гладкости $L$, но и величины относительной погрешности градиента $\alpha$ введем дополнительный параметр $\beta \in (0, 0.5]$ следующим образом: $\beta \coloneqq 0.5 - \alpha$. Предлагаемый алгоритм приведен ниже как алгоритм~\ref{alg:2}.
\begin{algorithm}
    \caption{Градиентный спуск с адаптивной настройкой $L$ и $\alpha$}
    \label{alg:2}
    \begin{algorithmic}
        \State
        \textbf{Вход:} $x^0$, $L_{min} \geq \mu > 0$, $L_0 \geq L_{min}$, $\alpha_{min} \in [0, 0.5)$, $\alpha_0 \in [\alpha_{min}, 0.5)$.
        \begin{enumerate}
            \renewcommand{\labelenumi}{\textbf{\arabic{enumi}:}}
            \item
                $k = 0$; $\beta_{max} = 0.5 - \alpha_{min}$; $\beta_0 = 0.5 - \alpha_0$.
            \item
                $L_{k+1} = \max\cbr*{\frac{L_k}{2}, L_{min}}$; $\beta_{k+1} = \min\cbr*{2\beta_k, \beta_{max}}$; $\alpha_{k+1} = 0.5 - \beta_{k+1}$.
            \item
                \begin{equation}\label{eq:alg2_step}
                    x^{k+1} = x^k - \frac{1}{L_{k+1}} \frac{1 - 2\alpha_{k+1}}{1 - \alpha_{k+1}} \Tilde{\nabla} f(x^k)
                \end{equation}
            \item Если
                \begin{multline}\label{eq:alg2_cond}
                    f(x^{k+1}) \leq f(x^k) + \abr*{\Tilde{\nabla} f(x^k), x^{k+1} - x^k} + \frac{L_{k+1}}{2} \norm*{x^{k+1} - x^k}^2 +\\+ \frac{\alpha_{k+1}}{1 - \alpha_{k+1}} \norm*{\Tilde{\nabla} f(x^k)} \norm*{x^{k+1} - x^k},
                \end{multline}
                то переходим на шаг 5. Иначе $L_{k+1} \coloneqq 2 L_{k+1}$, $\beta_{k+1} \coloneqq 0.5 \beta_{k+1}$, $\alpha_{k+1} \coloneqq 0.5 - \beta_{k+1}$ и возвращаемся на шаг 3.
            \item
                Если правило остановки не выполнено, то $k \coloneqq k+1$ и переходим на шаг 2.
        \end{enumerate}
        \textbf{Выход:} $x^k$.
    \end{algorithmic}
\end{algorithm}
Можно считать, что $L_0 \leq L$ и $\alpha_0 \leq \alpha$ (то есть $\beta_0 \geq \beta$), поэтому
\begin{equation}
    \max_{k<N} L_{k+1} \leq L \max\cbr*{2, \frac{\beta_{max}}{0.5 \beta}} = 2L \max\cbr*{1, \frac{0.5 - \alpha_{min}}{0.5 - \alpha}}
\end{equation}
и
\begin{equation}
    \min_{k<N} \beta_{k+1} \geq \frac{\beta}{\max\cbr*{2, \frac{2L}{L_{min}}}} = \frac{0.5 - \alpha}{2} \min\cbr*{1, \frac{L_{min}}{L}}.
\end{equation}
И, следовательно,
\begin{equation}
    \max_{k<N} \alpha_{k+1} = 0.5 - \min_{k<N} \beta_{k+1} \leq 0.5 - \frac{0.5 - \alpha}{2} \min\cbr*{1, \frac{L_{min}}{L}}.
\end{equation}
Введем соответствующие обозначения:
\begin{equation}\label{eq:L_max_2}
    L_{max} \coloneqq 2L \max\cbr*{1, \frac{0.5 - \alpha_{min}}{0.5 - \alpha}},
\end{equation}
\begin{equation}\label{eq:alpha_max}
    \alpha_{max} \coloneqq 0.5 - \frac{0.5 - \alpha}{2} \min\cbr*{1, \frac{L_{min}}{L}}.
\end{equation}

Проведем анализ получаемой алгоритмом \ref{alg:2} выходной точки в рамках следующего предположения, аналогичного используемому авторами \cite{kuruzov2023gradient} для случая абсолютной погрешности.
\begin{assumption}\label{ass:1}
    Алгоритм \ref{alg:2} может запрашивать значение градиента в любой текущей точке $x^k$ с произвольной относительной точностью, которая не выше чем $\alpha_{min}$, а конкретно с точностью $\alpha_{k+1}$.
\end{assumption}
\noindent Это позволяет записать PL-условие для неточного градиента в точке $x^k$, заменив $\alpha$ на $\alpha_{k+1}$:
\begin{equation}\label{eq:PL_rel_k}
    f(x^k) - f^* \leq \frac{1}{2 \mu (1 - \alpha_{k+1})^2} \norm*{
        \Tilde{\nabla} f(x^k)
    }^2.
\end{equation}
Условие \eqref{eq:PL_rel_k}, в свою очередь, позволяет, аналогично полученным в предыдущем подразделе результатам, записать следующую оценку скорости сходимости по функции для алгоритма \ref{alg:2}:
\begin{equation}\label{eq:final_2}
    f(x^N) - f^* \leq \prod\limits_{k = 0}^{N - 1} \rbr*{1 - \frac{\mu}{L_{k+1}} \xi_{k+1}} \rbr*{f(x^0) - f^*} \leq \rbr*{1 - \frac{\mu}{L_{max}} \xi_{max}}^N \rbr*{f(x^0) - f^*},
\end{equation}
где
\begin{equation}\label{eq:xi_2}
    \xi_{k+1} \coloneqq (1 - 2\alpha_{k+1})^2 \quad\text{и}\quad \xi_{max} \coloneqq (1 - 2\alpha_{max})^2,
\end{equation}
а также получить следующий основной результат этого раздела.
\begin{theorem}\label{th:2}
     Пусть функция $f$ удовлетворяет условию Липшица градиента \eqref{eq:L} и условию Поляка--Лоясиевича \eqref{eq:PL}. И пусть при работе алгоритма~\ref{alg:2} выполнилось одно из следующих двух условий:
     \begin{enumerate}
         \item алгоритм~\ref{alg:2} проработал $N_{**}$ итераций, где
             \begin{equation}\label{eq:N_2}
                 N_{**} = \ceil*{\frac{L_{max}}{\mu\xi_{max}} \log\rbr*{\frac{\mu \rbr*{f(x^0) - f^*}}{\varepsilon}}};
             \end{equation}
         \item для некоторого $N < N_{**}$ выполнено правило остановки в виде
            \begin{equation}
                \norm*{\Tilde{\nabla} f(x^N)}^2 \leq 2 \epsilon (1 - \alpha_{N+1})^2,
            \end{equation}
     \end{enumerate}
     где $\xi_{max}$ задается формулой \eqref{eq:xi_2}, а $L_{max}$ и $\alpha_{max}$ задаются формулами \eqref{eq:L_max_2} и \eqref{eq:alpha_max} соответственно. Тогда верны следующие оценки на выходную точку $\hat{x}$ алгоритма~\ref{alg:2} ($\hat{x} = x^{N_{**}}$ или $\hat{x} = x^{N}$ соответственно):
     \begin{equation}
         f(\hat{x}) - f^* \leq \frac{\epsilon}{\mu}
     \end{equation}
     и
     \begin{equation}\label{eq:traj_2}
         \norm*{\hat{x} - x^0} \leq \frac{2 L_{max}}{\mu\xi_{max}} \sqrt{\frac{2}{L_{min}} \rbr*{f(x^0) - f^*}}.
     \end{equation}
\end{theorem}

\begin{remark}
    Предположение~\ref{ass:1} довольно естественно для случая относительной погрешности, так как одной из основных областей её возникновения является безградиентная оптимизация, где зачастую можно регулировать точность вычисления приближенного значения градиента в запрашиваемой точке.
\end{remark}

\begin{remark}
В силу оценки \eqref{eq:traj_2} из теоремы \ref{th:2} для алгоритма \ref{alg:2} справедливы все выводы замечания \ref{Puch_rem1}, сделанного выше для алгоритма \ref{alg:1}. Это продемонстрировано ниже, в вычислительных экспериментах в разделе~\ref{sec:5}.
\end{remark}

\begin{remark}
    Исходя из предположения~\ref{ass:1}, алгоритм имеет доступ к градиенту с относительной погрешностью $\alpha_{min}$. Тем не менее понижение запрашиваемой точности до $\alpha_{k+1}$ ведет к ускорению вычисления неточного градиента, а также к возможному уменьшению второго настраиваемого параметра $L_{k+1}$, что имеет положительный эффект на оценку \eqref{eq:final_2} и оценки \eqref{eq:N_2}, \eqref{eq:traj_2} из теоремы \ref{th:2}.
\end{remark}

\begin{remark}
    Предположение~\ref{ass:1} не обязательно для получения линейной скорости сходимости. В случае когда алгоритм не может регулировать величину погрешности градиента, имеют место аналоги оценки \eqref{eq:final_2} и теоремы~\ref{th:2} для
    \begin{equation}
        \xi_{k+1} \coloneqq (1 - 2\alpha_{k+1})^2 \frac{(1 - \alpha)^2}{(1 - \alpha_{k+1})^2} \quad\text{и}\quad \xi_{max} \coloneqq (1 - 2\alpha_{max})^2 \frac{(1 - \alpha)^2}{(1 - \alpha_{max})^2}
    \end{equation}
    соответственно.
\end{remark}

\begin{remark}
    Во всех приведённых в настоящем подразделе статьи результатах $L_{max}$ и $\alpha_{max}$, оценивающие максимально возможные значения для $L_{k+1}$ и $\alpha_{k+1}$, могут быть заменены на $\max_{k < N} L_{k+1}$ и $\max_{k < N} \alpha_{k+1}$ соответственно.
\end{remark}

\begin{remark}
    Алгоритм~\ref{alg:2} (как и алгоритм~\ref{alg:1}) за счёт адаптивности применим даже для задач, не обладающих глобально свойством $L$-гладкости на всём пространстве. Это продемонстрировано ниже в вычислительных экспериментах для некоторых задач в разделе~\ref{sec:5}.
\end{remark}

\begin{remark}
    Можно оценить число повторений отдельно взятых шагов алгоритма~\ref{alg:2}. Если $L_{k+1} \geq L$ и $\alpha_{k+1} \geq \alpha$ (то есть $\beta_{k+1} \leq \beta$), то повторения 4-го шага заканчиваются. Поэтому повторений 4-го шага за все итерации не более $2N + \log_2 \rbr*{2 \max\cbr*{\frac{L}{L_{min}}, \frac{0.5 - \alpha_{min}}{0.5 - \alpha}}}$. Таким образом можно, например, оценить число обращений к оракулу на протяжении работы алгоритма.
\end{remark}

\section{Вычислительные эксперименты}\label{sec:5}

Для проведения экспериментального анализа работы итогового для настоящей статьи алгоритма \ref{alg:2} были выбраны функция Розенброка, а также функция Нестерова--Скокова. Обе эти функции невыпуклы и удовлетворяют PL-условию \eqref{eq:PL} на любом компакте (см., например, \cite{polyak2022stopping}). При этом для них нет возможности гарантировать условие $L$-гладкости на всём пространстве.

Все вычислительные эксперименты проводились на платформе Google Colab (\url{https://colab.research.google.com/}). Для генерации равномерного шума из многомерного шара использовался инструмент nengo.dists.UniformHypersphere из библиотеки nengo (\url{https://www.nengo.ai/}).

\subsection{Функция Розенброка}

\begin{figure}[h]
\centering
\includegraphics[width=\textwidth]{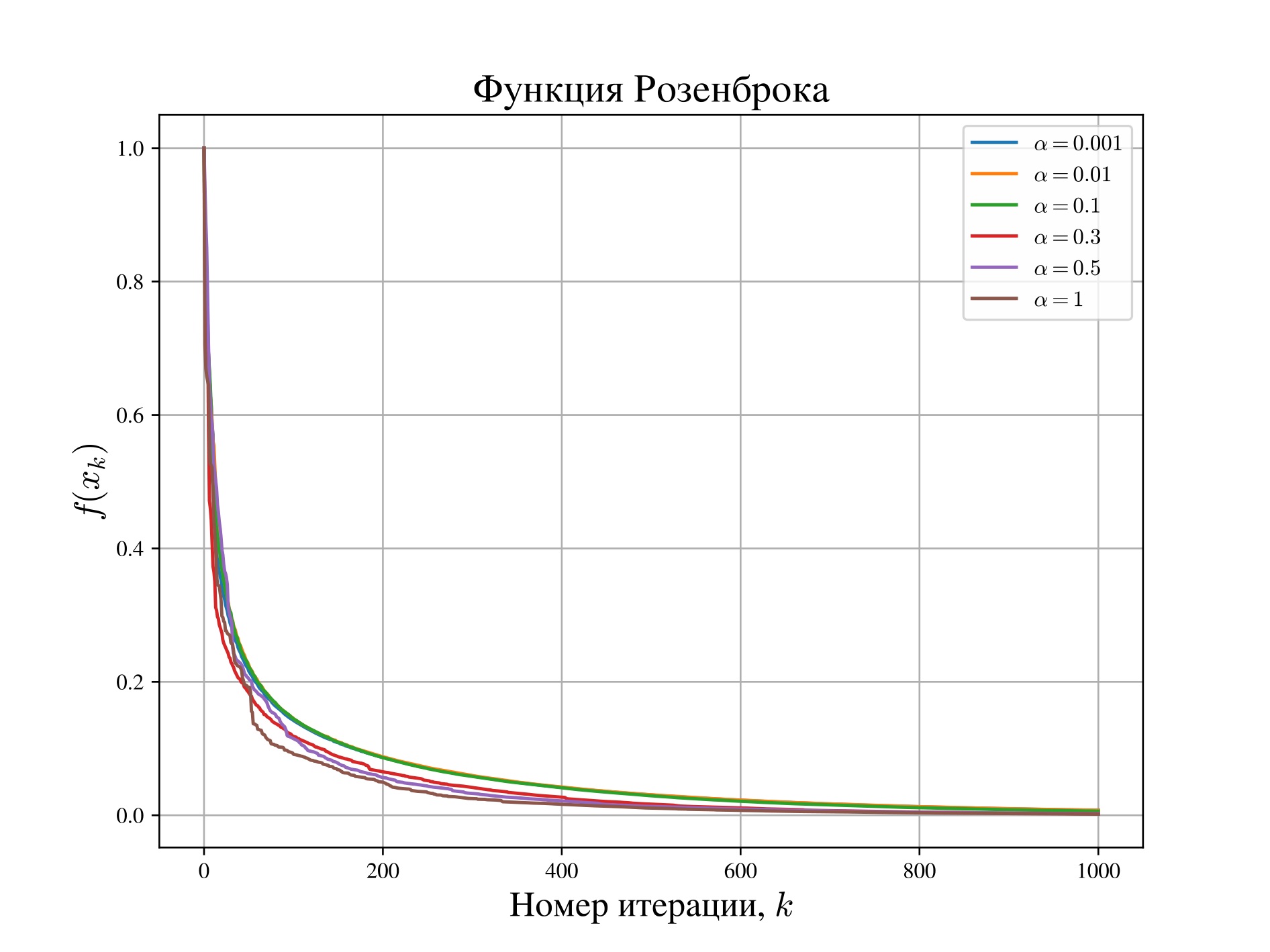}
\caption{График зависимости значения функции от номера итерации в зависимости от параметра помех $\alpha$ для функции Розенброка}
\label{fig:1}
\end{figure}

Функцией Розенброка называется следующая функция двух переменных
\begin{equation}
    f(x_1, x_2) = 100 \rbr*{x_2 - x_1^2}^2 + (x_1 - 1)^2.
\end{equation}
Как легко видеть, ее глобальный минимум находится в точке $(x_1, x_2) = (1, 1)$, и $f^* = 0$. На вход алгоритму подавались следующие входные данные: $x^0 = (0, 0)$, $L_{min} = 0.01$, $L_0 = 1$, $\alpha_{min} = 0.001$, $\alpha_0 = 0.01$. В качестве неточного градиента рассматривалось значение градиента $\nabla f(x_1, x_2)$, зашумленное равномерно распределенным в $2$-мерном шаре радиуса $\alpha \lVert \nabla f(x_1, x_2) \rVert$ относительным шумом, где параметр $\alpha$ варьировался. Результаты эксперимента приведены на рис.~\ref{fig:1}, а также в табл.~\ref{tab:1} и \ref{tab:2}.

\begin{table}[ht!]
    \caption{Результат работы алгоритма \ref{alg:2} после $N = 1000$ итераций для функции Розенброка\centering}
    \label{tab:1}
    \centering
    \begin{tabular}{c||c|c|c|c|c|c}
        $\alpha$ & $0.001$ & $0.01$ & $0.1$ & $0.3$ & $0.5$ & $1$\\ \hline
        $f(x_N)$ & $0.0074$ & $0.0075$ & $0.0060$ & $0.0021$ & $0.0018$ & $0.0017$
    \end{tabular}
\end{table}


\begin{table}[ht!]
    \caption{Результат работы алгоритма \ref{alg:2} после $N = 10000$ итераций для функции Розенброка\centering}
    \label{tab:2}
    \centering
    \begin{tabular}{c||c|c|c|c|c|c}
        $\alpha$ & $0.001$ & $0.01$ & $0.1$ & $0.3$ & $0.5$ & $1$ \\ \hline
        $f(x_N)$ & $1.5 \times 10^{-19}$ & $1.3 \times 10^{-19}$ & $1.6 \times 10^{-19}$ & $2.6 \times 10^{-16}$ & $2.7 \times 10^{-15}$ & $7.3 \times 10^{-17}$
    \end{tabular}
\end{table}

Для всех рассмотренных $\alpha$ метод уже к $10000$-ой итерации сходится по функции и, следовательно, по аргументу к глобальному минимуму с машинной точностью. Отметим, что предполагаемое замедление сходимости при увеличении $\alpha$ не наблюдается как минимум до $1000$-ой итерации (см.~табл.~\ref{tab:1}), однако по итогу при сильном приближении к минимуму это начинает незначительно сказываться, что видно из табл.~\ref{tab:2}. В том числе, оказалось приемлемым значение $\alpha = 1$, которое в теоретических рассуждениях нами отбрасывалось.

\subsection{Функция Нестерова--Скокова}

Функцией Нестерова--Скокова называется следующая функция $n$ переменных
\begin{equation}
    f(x_1, x_2, \ldots, x_n) = \frac{1}{4} (1 - x_1)^2 + \sum\limits_{i = 1}^{n - 1} \rbr*{x_{i+1} - 2 x_i^2 + 1}^2.
\end{equation}
Эту функцию также называют обобщением функции Розенброка на многомерный случай. Как легко видеть, ее глобальный минимум находится в точке $(x_1, x_2, \ldots, x_n) = (1, 1, \ldots, 1)$, и $f^* = 0$. На вход алгоритму подавались следующие входные данные: $L_{min} = 0.01$, $\alpha_{min} = 0.001$, $\alpha_0 = 0.01$. В качестве неточного градиента рассматривалось значение градиента $\nabla f(x_1, x_2, \ldots, x_n)$, зашумленное равномерно распределённым в $n$-мерном шаре радиуса $\alpha \lVert \nabla f(x_1, x_2, \ldots, x_n) \rVert$ относительным шумом, где параметр $\alpha$ варьировался. Также варьировались $x^0$ и $L_0$. Размерность $n$ выбрана равной $100$.

\begin{figure}[h]
\centering
\includegraphics[width=\textwidth]{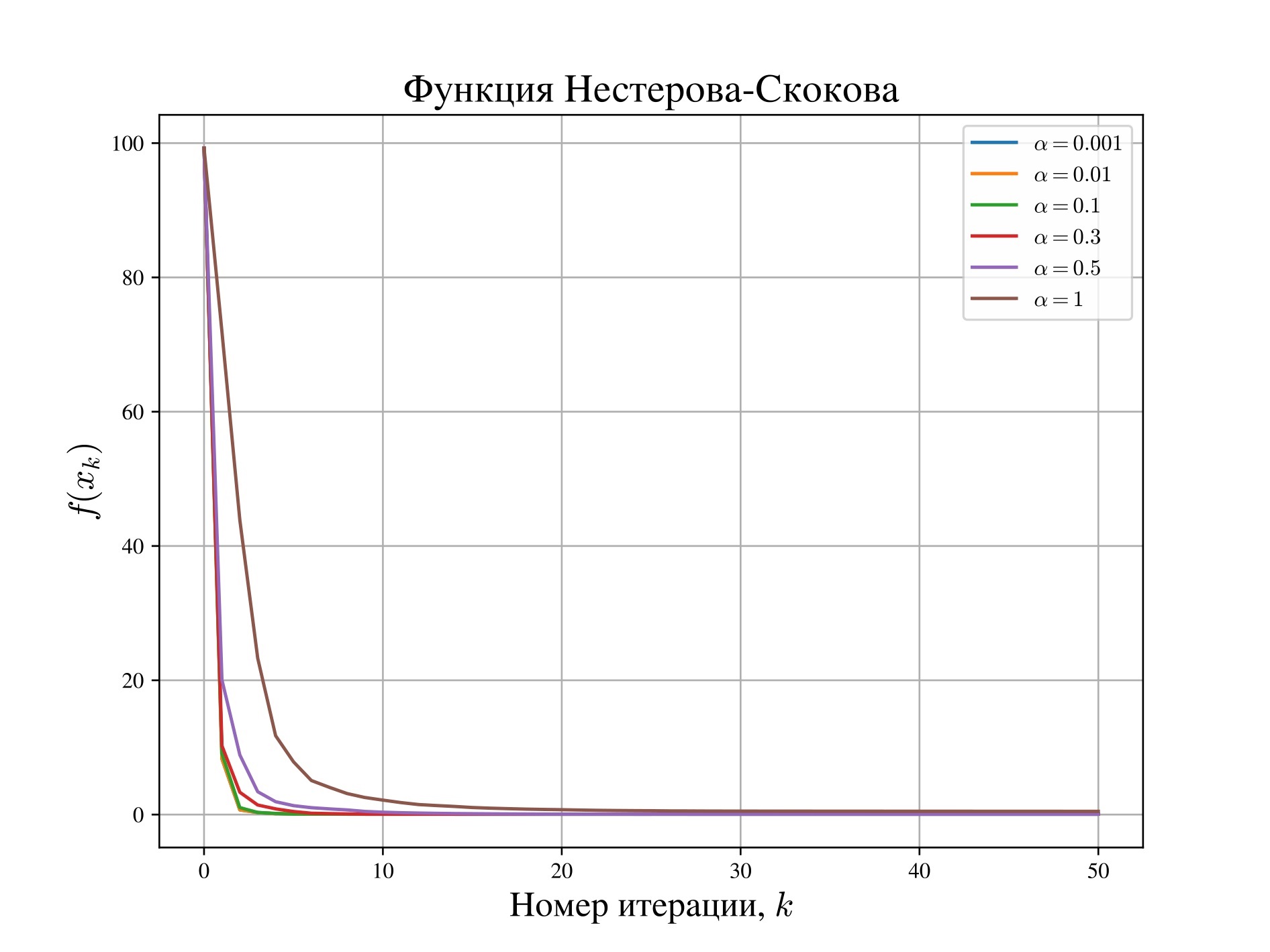}
\caption{График зависимости значения функции от номера итерации в зависимости от параметра помех $\alpha$ для функции Нестерова--Скокова ($x^0 = (0, 0, \ldots, 0)$, $L_0 = 1$)}
\label{fig:2}
\end{figure}

Результаты эксперимента приведены на рис.~\ref{fig:2}, а также в табл.~\ref{tab:3}, \ref{tab:4}, \ref{tab:5} и \ref{tab:6}.

\begin{table}[ht!]
    \caption{Результат работы алгоритма \ref{alg:2} после $N = 10$ итераций для функции Нестерова--Скокова ($x^0 = (0, 0, \ldots, 0)$, $L_0 = 1$)\centering}
    \label{tab:3}
    \centering
    \begin{tabular}{c||c|c|c|c|c|c}
        $\alpha$ & $0.001$ & $0.01$ & $0.1$ & $0.3$ & $0.5$ & $1$\\ \hline
        $f(x_N)$ & $0.058$ & $0.058$ & $0.059$ & $0.073$ & $0.261$ & $2.631$
    \end{tabular}
\end{table}

\begin{table}[ht!]
    \caption{Результат работы алгоритма \ref{alg:2} после $N = 50$ итераций для функции Нестерова--Скокова ($x^0 = (0, 0, \ldots, 0)$, $L_0 = 1$)\centering}
    \label{tab:4}
    \centering
    \begin{tabular}{c||c|c|c|c|c|c}
        $\alpha$ & $0.001$ & $0.01$ & $0.1$ & $0.3$ & $0.5$ & $1$\\ \hline
        $f(x_N)$ & $0.058$ & $0.058$ & $0.058$ & $0.058$ & $0.058$ & $0.058$
    \end{tabular}
\end{table}

\begin{table}[ht!]
    \caption{Результат работы алгоритма \ref{alg:2} после $N = 10$ итераций для функции Нестерова--Скокова ($x^0 = (-1, 1, \ldots, 1)$, $L_0 = 0.1$)\centering}
    \label{tab:5}
    \centering
    \begin{tabular}{c||c|c|c|c|c|c}
        $\alpha$ & $0.001$ & $0.01$ & $0.1$ & $0.3$ & $0.5$ & $1$\\ \hline
        $f(x_N)$ & $1.2 \times 10^{-6}$ & $6.7 \times 10^{-5}$ & $0.98$ & $0.98$ & $0.98$ & $0.98$
    \end{tabular}
\end{table}

\begin{table}[ht!]
    \caption{Результат работы алгоритма \ref{alg:2} после $N = 50$ итераций для функции Нестерова--Скокова ($x^0 = (-1, 1, \ldots, 1)$, $L_0 = 0.1$)\centering}
    \label{tab:6}
    \centering
    \begin{tabular}{c||c|c|c|c|c|c}
        $\alpha$ & $0.001$ & $0.01$ & $0.1$ & $0.3$ & $0.5$ & $1$\\ \hline
        $f(x_N)$& $4.4 \times 10^{-11}$ & $3.2 \times 10^{-9}$ & $0.98$ & $0.98$ & $0.98$ & $0.98$
    \end{tabular}
\end{table}

Для случая $x^0 = (0, 0, \ldots, 0)$, $L_0 = 1$ метод для всех рассмотренных $\alpha$ быстро сходится к минимуму, но не глобальному, а только локальному, который, однако, относительно близок к глобальному по функции. Для случая $x^0 = (-1, 1, \ldots, 1)$, $L_0 = 0.1$ при $\alpha = 0.001$ и $\alpha = 0.01$ метод быстро сходится уже к глобальному минимуму, однако при остальных рассмотренных $\alpha$ метод почти сразу сваливается в локальный минимум с далеким от оптимального значением функции. В обоих случаях, в отличие от функции Розенброка, здесь отчётливо заметно снижение скорости сходимости при увеличении параметра $\alpha$, пусть и не значительное.

\section{Заключение} В статье продолжены начатые недавно в работах \cite{polyak2022stopping, kuruzov2023gradient} исследования по анализу поведения траекторий адаптивных методов градиентного типа для задач с условием Поляка--Лоясиевича в предположении доступности методу в текущей точке информации о целевой функции. Интерес к этому классу задач обусловлен тем, что к нему сводятся важные вопросы исследования нелинейных систем с перепараметризацией в глубоком обучении \cite{belkin2021fit} и многие другие прикладные задачи \cite{karimi2016linear}. Такие задачи могут быть невыпуклыми \cite{polyak2020new}, но при этом можно гарантировать для них сходимость методов градиентного типа со скоростью геометрической прогрессии в случае $L$-гладкости. Если в статьях \cite{polyak2022stopping, kuruzov2023gradient} исследовался случай аддитивной неточности градиента и аддитивной неточности задания целевой функции, то в настоящей работе рассмотрена ситуация доступности в текущей точке градиента с относительной погрешностью, что актуально, например, с точки зрения безградиентных методов \cite{berahas2022theoretical}. Принципиальное отличие полученных результатов от \cite{polyak2022stopping, kuruzov2023gradient}~--- доказательство сохранения сходимости предложенных вариаций градиентного метода со скоростью геометрической прогрессии при некоторых предположениях об уровне относительной погрешности. В случае абсолютной неточности информации о градиенте в \cite{polyak2022stopping, kuruzov2023gradient} удалось лишь оценить (показана оптимальность таких оценок) степень отклонения от сходимости со скоростью геометрической прогрессии в случае известной оценки параметра погрешности информации о градиенте и предложить правила ранней остановки методов, гарантирующие достижение оптимального уровня качества приближённого решения по функции. При этом адаптивность рассмотренных методов открывает возможности их применения в задачах, для которых нет глобального свойства $L$-гладкости, в том числе и в негладких по аналогии с методикой универсальных градиентных методов Ю.Е.~Нестерова \cite{nesterov2014universal}.

\hfil \hbox to 0.3\textwidth{\hrulefill} \smallskip

Благодарности. Работа выполнена при поддержке гранта Российского научного фонда и города Москвы № 22-21-20065
(\url{https://rscf.ru/project/22-21-20065/}).

\end{document}

%% file: header.tex
\usepackage[utf8]{inputenc} 
\usepackage[T2A]{fontenc}    
\usepackage[english,russian]{babel}
\usepackage{mmap}
\usepackage{hyperref}       
\usepackage{url}            
\usepackage{booktabs}       
\usepackage{amsfonts}       
\usepackage{nicefrac}       
\usepackage{lipsum}
\usepackage{mathtools}
\usepackage{amsmath,amssymb,amsthm}
\usepackage{algorithm}
\usepackage{algpseudocode}
\usepackage{pifont}
\usepackage{cases}
\usepackage{graphicx}
\usepackage{stackengine}    
\usepackage{titlesec}
\usepackage[none]{hyphenat}
\usepackage{subfigure}
\usepackage[labelsep=period]{caption}
\hypersetup{
  colorlinks   = true,    
  urlcolor     = green,   
  linkcolor    = blue,    
  citecolor    = red      
}

\def\keywordname{{\bfseries Ключевые слова:}}%
\def\keywords#1{\par\addvspace\medskipamount{\rightskip=0pt plus1cm
\def\and{\ifhmode\unskip\nobreak\fi\ $\cdot$ }\noindent\keywordname\enspace\ignorespaces#1\par}}

\def\keywordnameeng{{\bfseries Keywords:}}%
\def\keywordseng#1{\par\addvspace\medskipamount{\rightskip=0pt plus1cm
\def\and{\ifhmode\unskip\nobreak\fi\ $\cdot$ }\noindent\keywordnameeng\enspace\ignorespaces#1\par}}

\titlelabel{\thetitle. }

%% file: main_ArXiv.bbl
\begin{thebibliography}{99}

\bibitem{beck2017first} \textit{Beck A.} First-Order Methods in Optimization.~--- Society for Industrial and Applied Mathematics (Oct 2017). https://doi.org/10.1137/1.9781611974997

\bibitem{nesterov2004introductory} \textit{Nesterov Y.} Introductory Lectures on Convex Optimization.~--- Springer US (2004). https://doi.org/10.1007/978-1-4419-8853-9

\bibitem{polyak1963gradient} \textit{Polyak B.T.} Gradient methods for the minimisation of functionals~// USSR Computational Mathematics and Mathematical Physics 3(4), 864--878 (Jan 1963). https://doi.org/10.1016/0041-5553(63)90382-3

\bibitem{Yue2022low} \textit{Yue P., Fang C., Lin Z.} On the lower bound of minimizing polyak- lojasiewicz functions. // The Thirty Sixth Annual Conference on Learning Theory. --- Proceedings of Machine Learning Research 195, 2948 -- 2968. (2023). https://proceedings.mlr.press/v195/yue23a/yue23a.pdf

\bibitem{devolder2013first} \textit{Devolder O., Glineur F., Nesterov Y.} First-order methods of smooth convex optimization with inexact oracle~// Mathematical Programming 146(1-2), 37--75 (Jun 2013). https://doi.org/10.1007/s10107-013-0677-5

\bibitem{polyak1983introduction} \textit{Поляк Б.Т.} Введение в оптимизацию.~--- М.: Наука, 1983.

\bibitem{polyak2022stopping} \textit{Polyak B.T., Kuruzov I.A., Stonyakin F.S.} Stopping rules for gradient methods for non-convex problems with additive noise in gradient~// Journal of Optimization Theory and Applications 198, 531–551 (2023). https://doi.org/10.1007/s10957-023-02245-w

\bibitem{kuruzov2023gradient} \textit{Kuruzov I.A., Stonyakin F.S., Alkousa M.S.} Gradient-type methods for optimization problems with Polyak--{\L}ojasiewicz condition: Early stopping and adaptivity to inexactness parameter~// Olenev, N., Evtushenko, Y., Jaćimović, M., Khachay, M., Malkova, V., Pospelov, I. (eds) Advances in Optimization and Applications. OPTIMA 2022. Communications in Computer and Information Science 1739, 18 -- 32(2022). https://doi.org/10.48550/ARXIV.2212.04226

\bibitem{gasnikov2021mccme} \textit{Гасников А.В.} Современные численные методы оптимизации. Метод универсального градиентного спуска.~--- М.: МЦНМО, 2021.

\bibitem{carter1991global} \textit{Carter R.G.} On the global convergence of trust region algorithms using inexact gradient information~// SIAM Journal on Numerical Analysis 28(1), 251--265 (Feb 1991). https://doi.org/10.1137/0728014

\bibitem{berahas2022theoretical} \textit{Berahas A.S., Cao L., Choromanski K., Scheinberg K.} A theoretical and empirical comparison of gradient approximations in derivative-free optimization~// Foundations of Computational Mathematics 22(2), 507--560 (May 2021). https://doi.org/10.1007/s10208-021-09513-z

\bibitem{belkin2021fit} \textit{Belkin M.} Fit without fear: remarkable mathematical phenomena of deep learning through the prism of interpolation~// Acta Numerica 30, 203--248 (2021). https://doi.org/10.1017/S0962492921000039

\bibitem{karimi2016linear} \textit{Karimi H., Nutini J., Schmidt M.} Linear convergence of gradient and proximal-gradient methods under the polyak-{\l}ojasiewicz condition~// Springer, 795--811 (2016). https://doi.org/10.1007/978-3-319-46128-1\_50

\bibitem{polyak2020new} \textit{Polyak B., Tremba A.} New versions of Newton method: step-size choice, convergence domain and under-determined equations~// Optimization Methods and Software 35(6), 1272--1303 (2020). https://doi.org/10.1080/10556788.2019.1669154

\bibitem{nesterov2014universal} \textit{Nesterov Y.} Universal gradient methods for convex optimization problems~// Mathematical Programming 152(1--2), 381--404 (May 2014). https://doi.org/10.1007/s10107-014-0790-0

\end{thebibliography}
